\documentclass{article}

\usepackage{latexsym}
\usepackage{amssymb,amsmath,amsopn}
\usepackage[dvips]{graphicx}
\usepackage{color,epsfig}
\usepackage{url}
\usepackage{multicol}
\usepackage{amsthm}

\pdfpagewidth \paperwidth
\pdfpageheight \paperheight
\newtheorem*{remark}{Remark}

\begin{document}
\title{\textbf{{Triangular fractal approximating \\ graphs and their covering \\paths and cycles}}}

\author{Kaszanyitzky, Andr\'as \\ kaszi75@gmail.com}
\date{}
\maketitle

\textbf{Abstract.} We observed and described the \emph{generalized Sierpi\'{n}ski Arrowhead Curve} in our previous paper [K17a]. Now we focus on its background structure. 
In \emph{Section 1} we summarize our previous results on the triangular grid and supplement them with Hamiltonian-cycles, tiling-cycles and a new kind of path on the possible largest trapezoid grid which are needed for the following sections. We describe the basic rule of the transformability of the paths and the cycles into each other and extend our grids to larger graphs. In \emph{Sections 2} and \emph{3} we define two kinds of graphs related to a checked fractal pattern on the \emph{generalized Sierpi\'{n}ski Gasket}. We continue our observations with the basic properties of these triangular fractal approximating graphs independently of the recursive curves. We will describe the numbers of their vertices and edges, and their covering paths and cycles in general case with recursive and explicit formulas. Some of their cardinality specify new integer sequences. We also find the bijective relations between these formations.

\section{Paths and cycles on the triangular grid}

In this section we summarize and supplement our previous definitions [K17a] and show a table with our results, the cardinality of these formations on simple triangular grids. We describe the basic rules of the transformability of the paths and cycles into each other and extend our grids to larger graphs.

\subsection{Checked generator pattern}

First we make a checked pattern on a \emph{triangular grid} of order $n$ by colouring the tiles that face upwards dark and colouring the rest of the subtriangles white. Then we substitute all the dark tiles with the contracted copy of this generator pattern. Our pattern is related to the two-dimensional \emph{generalized Sierpi\'{n}ski Gasket} $SG_{2,n}(k)$. We will be referring to them as $F_n$ generator pattern and as $F_n(k)$ fractal approximating pattern, if $k>1$.

\subsection{Triangular grids and their paths and cycles}

Our generator pattern $F_n$ contains $n^2$ subtriangles, and $T_{n-1}$ white tiles, $T_n$ dark tiles and $T_{n+1}$ grid points as three consecutive triangular numbers. The centroids of the dark tiles form the \emph{inscribed grid}, and their corners form the \emph{overall grid}. All of our paths originate from the leftmost node and terminate in the rightmost node of these grids.

Let us consider a self-avoiding tiling-path called \emph{S-path} (referring to Sierpi\'{n}ski), denoted by $S_n$, and a self-avoiding tiling-cycle called \emph{D-cycle}, denoted by $D_n$, on the overall grid. Both consist of $T_n$ edges. All of the edges must be lying on different dark subtriangles. For practical reasons we will be using the notation of McKenna: marking the tiles with little ticks in the middle of the edges [McK94]. See the left side of \emph{Figure 1}.

We denote the Hamiltonian-paths (H-paths) by $H_n$, and the Hamiltonian-cycles (C-cycles) by $C_n$ on the inscribed grid. They have a subset in which all edges have well-formed turns. We will describe this well-formed property later. These paths and cycles are bijective pairs and they have the same cardinality. They are unambiguously transformable into each other [K17a]. We call the well-formed Hamiltonian-paths W-path, and denote them by $W_n$. You can see a well-formed Hamiltonian-cycle on the middle of \emph{Figure 1}.

The so-called Z-paths $(Z_n)$  on the possible largest trapezoid grid without the uppermost node of the inscribed grid are also needed for \emph{Section 3}. See the right side of \emph{Figure 1}.

We enumerated these cardinality with our computer program, which is a smart backtrack algorithm. H-paths and Z-paths appear in [SEH05], W-paths first appear in [K17a], C-cycles appear in [P14] and [OEIS1], and D-cycles, which specify a new integer sequence, first appear here. See \emph{Table 1}.

\[
\begin{tabular}{|c||c||c||c||c||c||c|}
\hline
\bf
{\small n} & \bf {\small $T_n$} & \bf {\small $H_n$} & \bf {\small $W_n=S_n$ } & \bf {\small $Z_n$} & \bf {\small $C_n$} & \bf {\small $D_n$} \\ \hline \hline
\bf 2 & 3 & 1 &  1 & 1 & 1 & 1 \\ \hline
\bf 3 & 6 & 2 & 2 & 3 & 1 & 1 \\ \hline
\bf 4 & 10 & 10 & 4 & 11 & 3 & 3 \\ \hline
\bf 5 & 15 & 92 & 16 & 112 & 26 & 8 \\ \hline
\bf 6 & 21 & 1852 &  68 & 2286 & 474 & 42 \\ \hline
\bf 7 & 28 & 78032 &  464 & 94696 & 17214 & 240 \\ \hline
\bf 8 & 36 & 6846876 &  3828 & 8320626 & 1371454 & 2120 \\ \hline
\bf 9 & 45 & 1255156712 &  44488 & 1527633172 & 231924780 & 22724 \\ \hline
\end{tabular}
\]

\smallskip

\centerline{\bf \emph{Table 1.} Cardinality of H-, W-, S-, Z-paths and C- and D-cycles}
\centerline{\bf on the generator pattern $F_n$ consisting of $T_n$ dark tiles.}

\bigskip

\begin{figure}[ht]
\centering
\includegraphics[width=0.8\textwidth]{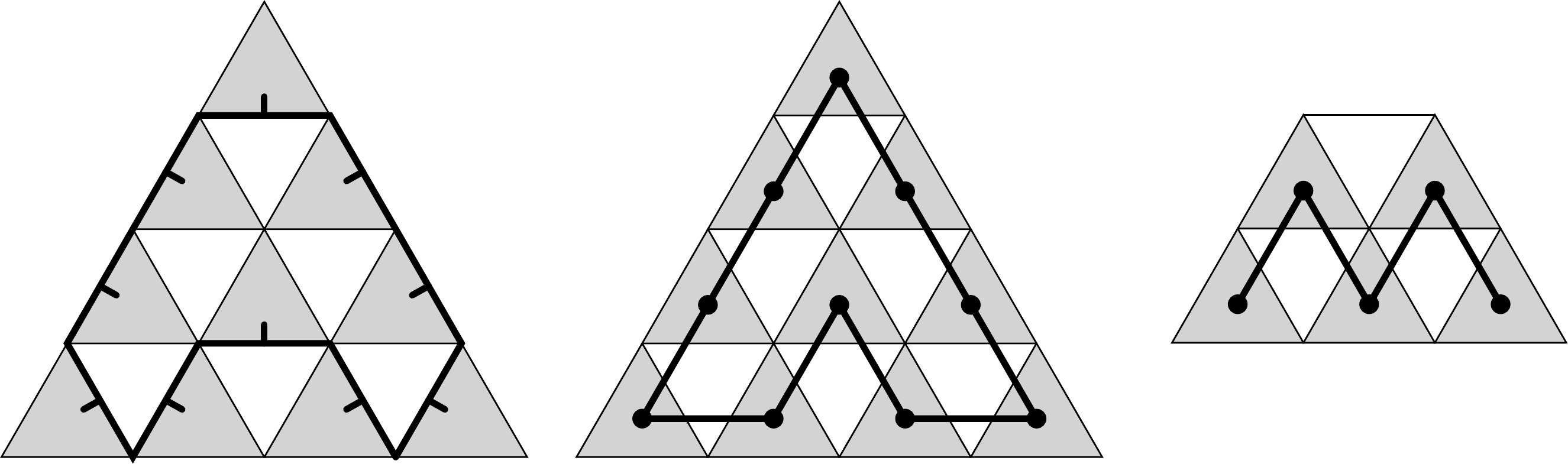}
\centerline{}
\centerline{\bf Figure 1. A tiling-cycle $(D_4)$ (left side), the corresponding} 
\centerline{\bf well-formed Hamiltonian-cycle $(C_4)$ (middle),}
\centerline{\bf and a Z-path $(Z_3)$ with a wrong turn of the edges (right side).}
\end{figure}

\subsection{Untransformable wrong turns of the edges}

Let us consider all edges of the paths and cycles on the inscribed grid described by a string of their absolute direction code consisting of $0$ to $5$ values in counterclockwise from the right. The direction right $ \pm 120^\circ$ means an even value and the direction left $ \pm 120^\circ$ means an odd value of the string.

If a path or a cycle on the inscribed grid consists of only well-formed turns, then it is unambiguously transformable into a tiling-path or a tiling-cycle.

The forbidden turns as $(d_i, d_{i+1})$ number pairs are the following: 

\begin{equation*}
d_{i+1}\not\equiv
\begin{cases}
(d_i +4) \mod 6 & \qquad \text{if $d_i$ is even} \\
(d_i +2) \mod 6 & \qquad \text{if $d_i$ is odd}
\end{cases}
\end{equation*}

which means that the next edge cannot turn $120^{\circ}$ to the right after an even direction and it cannot turn $120^{\circ}$ to the left after an odd direction. Naturally, turning back by $180^{\circ}$ is also forbidden.

For example on the right side of \emph{Figure 1}, Z-path contains a wrong turn of the edges (the middle edge-pair), therefore it is an untransformable path. These three dark tiles have only one contact point instead of two, therefore we cannot connect them with three consecutive edges of a tiling-path.

\subsection{Extending our grids to larger graphs}

The \emph{generalized Sierpi\'{n}ski Gasket} fractal family contains two kinds of \emph{triangular fractal approximating graphs} as the background structure of our recursive curves. They are the extended version of the \emph{inscribed grid} and the \emph{overall grid} in larger approximations, where $k>1$.

We observe triangular graphs based on the $k$-th power of the $n$-th \emph{Triangular number}, denoted by $T_n^k$. By connecting the centroids of the neighbouring dark tiles of $F_n(k)$ we get a graph that we call the \emph{Inscribed Graph} denoted by $I_n^k$. By connecting the corners of the neighbouring dark tiles of $F_n(k)$ we get a graph that we call the \emph{Overall Graph} denoted by $O_n(k)$. We will observe and describe their properties in the rest of this paper.

\section{The Overall Graph $O_n(k)$}

Let us define the \emph{Overall Graph $(O_n(k))$}, related to the $F_n(k)$ fractal approximating pattern where we replace all the dark tiles with their corners as the nodes of the graph and with their sides as the edges of the graph. 

It consists of $T_n^{k-1}$ simple $T_{n+1}$ sized triangular grids which share their corners with their neighbour grids. See \emph{Figure 2}. We will describe the numbers of their nodes and edges, and the numbers of their paths and cycles in this section.

$O_2(k)$ with all its possible connecting edges is also known as the \emph{Sierpi\'{n}ski Sieve Graph}.

\begin{figure}[ht]
\centering
\includegraphics[width=0.6\textwidth]{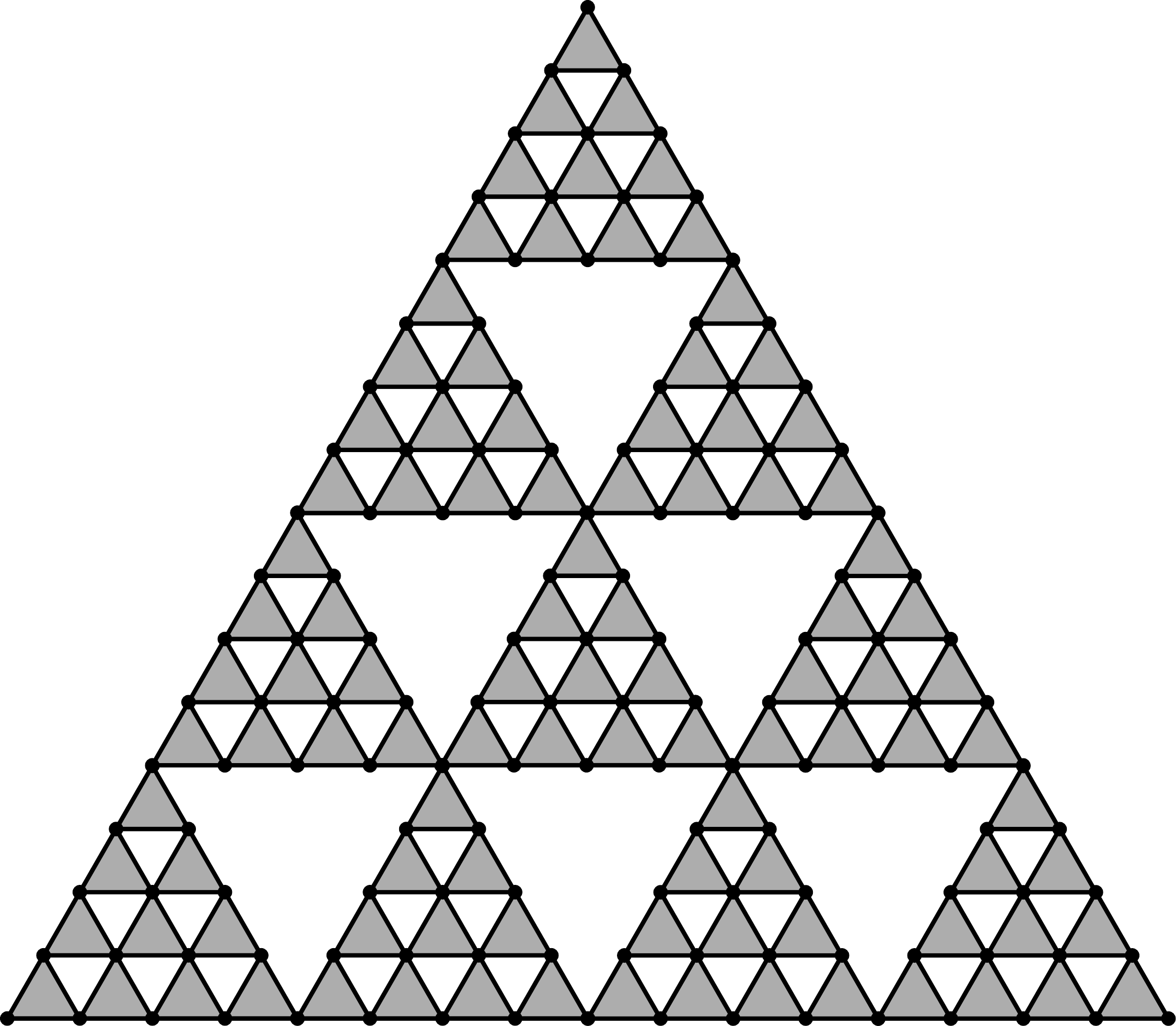}
\centerline{}
\centerline{\bf Figure 2. The Overall Graph $O_4(2)$}
\centerline{has $\left\vert O_4(2) \right\vert = 135$ nodes, $E(O_4(2))=300$ edges} \centerline{and $S_{4,2}=4^{11}=4194304$ possible S-paths.}
\end{figure}

\subsection{The numbers of the nodes and the edges in $O_n(k)$}

The overall grid consists of $T_{n+1}$ grid points, therefore  $\left\vert O_n(1) \right\vert=T_{n+1}$. 

In further approximations we substitute the dark tiles with $T_n$ smaller overall grids which share their corners. By summerizing the nodes of the $T_n$ smaller grids we counted the common nodes twice on each side of the overall graph and we counted the common nodes 3 times inside the overall graph, so we have to subtract these values from the result:

\begin{center}
$\left\vert O_n(k) \right\vert \quad = \quad {\left\vert O_n(k-1) \right\vert \cdot T_n} - 3(n-1) -2T_{n-2} \quad = \quad {\left\vert O_n(k-1) \right\vert \cdot T_n} - n^2 + 1$
\end{center}

We can transform this recursive formula to explicit formula and we get integer sequences for each $n$. 

\begin{center}
$\left\vert O_n(k) \right\vert = \dfrac{(n+4) \bigg( \dfrac{n(n+1)}{2} \bigg)^k +2(n+1)}{n+2}$
\end{center}

We always get our result in this form: $\left\vert O_n(k) \right\vert = \dfrac{aT_n^k+b}{c}$ where $a$, $b$ and $c$ values can be simplified by 2 for each even values of $n$.

\bigskip

For example: $\left\vert O_3(k) \right\vert=\dfrac{7 \cdot 6^k + 8}{5}$, \qquad \quad $\left\vert O_4(k) \right\vert=\dfrac{4 \cdot 10^k + 5}{3}$, 

\bigskip

\qquad \qquad \qquad $\left\vert O_5(k) \right\vert=\dfrac{9 \cdot 15^k + 12}{7}$,
\qquad $\left\vert O_6(k) \right\vert=\dfrac{5 \cdot 21^k + 7}{4}$, \qquad etc.

\bigskip

See the first 6 values of these integer sequences in \emph{Table 2}.

\[
\begin{tabular}{|c||c|c|c|c|c|c|}
\hline
\bf
{\small $\left\vert O_n(k) \right\vert$} & \bf { \small $k=1$} & \bf {\small $2$} & \bf {\small $3$ } & \bf {\small $4$} & \bf {\small $5$} & \bf {\small $6$}\\ \hline \hline
\bf $n=2$ & 6 & 15 &  42 & 123 & 366 & 1095 \\ \hline
\bf 3 & 10 & 52 & 304 & 1816 & 10888 & 65320 \\ \hline
\bf 4 & 15 & 135 & 1335 & 13335 & 133335 & 1333335 \\ \hline
\bf 5 & 21 & 291 & 4341 & 65091 & 976341 & 14645091 \\ \hline
\bf 6 & 28 & 553 &  11578 & 243103 & 5105128 & 107207653 \\ \hline
\end{tabular}
\]

\smallskip

\centerline{\bf \emph{Table 2.} $\left\vert O_n(k) \right\vert = $ the number of the nodes in the overall graph.}

\begin{remark}
First row of \emph{Table 2} is known as sequence $A067771$ [OEIS2]. Second and third rows are also known [CC06]. Our explicit formula gives new integer sequences for $\left\vert O_n(k) \right\vert$, where $n>4$.
\end{remark}

The number of the edges in the overall graph is: \qquad $E(O_n(k))=3{T_n^k}$.

\subsection{S-paths on $O_n(k)$}

We denote S-paths on the overall graph by $S_{n,k}$. In \emph{Section 1} we enumerated $S_n$ values on the overall grid therefore $S_{n,1}=S_n$. See \emph{Table 1}.

In the second approximation we can use all the $S_n$ paths at $T_n$ places and we have $S_n$ ways to connect them, therefore $S_{n,2}=S_n^{T_n+1}$. 

\bigskip

\centerline{In general case for $k>1$: \qquad $S_{n,k}=S_n^{{(T_n+1)}^{k-1}}$}

\bigskip

On $O_2(k)$ we get a unique path $(S_{2,k}=1)$ for all $k$ values. This trivial case is the $k$-th approximation of the edge-rewriting \emph{Sierpi\'{n}ski Arrowhead Curve}.

For $n>2$ and $k>1$ we get: 
\smallskip

$S_{3,k}=2^{7^{k-1}} \qquad S_{4,k}=4^{11^{k-1}} \qquad S_{5,k}=16^{16^{k-1}} \qquad S_{6,k}=68^{22^{k-1}}$

\subsection{Tiling-cycles on $O_n(k)$}

To find the number of the tiling-cycles on the overall graph we have to substitute the connection of the smaller grids with tiling-cycles $(D_n)$ instead of S-paths. By modifying our previous formula we get the following.

The number of the tiling-cycles on $O_n(k)$ in general case with recursion is: 
\smallskip

\centerline{$D_{n,k}=D_n \cdot S_{n,k-1}^{T_n^{k-1}}$}

\bigskip

\centerline{The explicit formula is: \qquad $D_{n,k}=D_n \cdot \bigg( S_n^{{(T_n+1)}^{k-2}} \bigg)^{T_n^{k-1}}$}

\bigskip

\centerline{We get very large numbers:}

\renewcommand{\arraystretch}{1.4}
\[
\begin{tabular}{|c||c|c|c|}
\hline
\bf
{\small $D_{n,k}$} & \bf { \small $k=1$} & \bf {\small $2$} & \bf {\small $3$ } \\ \hline \hline
\bf $n=2$ & 1 & 1 & 1 \\ \hline
\bf 3 & 1 & $2^6$ & $2^{252}$ \\ \hline
\bf 4 & 3 & $3 \cdot 4^{10}$ & $3 \cdot 4^{1100}$  \\ \hline
\bf 5 & 8 & $8 \cdot 16^{15}$ & $2^{14403}$ \\ \hline
\end{tabular}
\]

\bigskip

\centerline{\bf \emph{Table 3.} The number of the tiling-cycles =}
\centerline{\bf $D_{n,k}$ values on the overall graph.}

\begin{figure}[ht]
\centering
\includegraphics[width=0.9\textwidth]{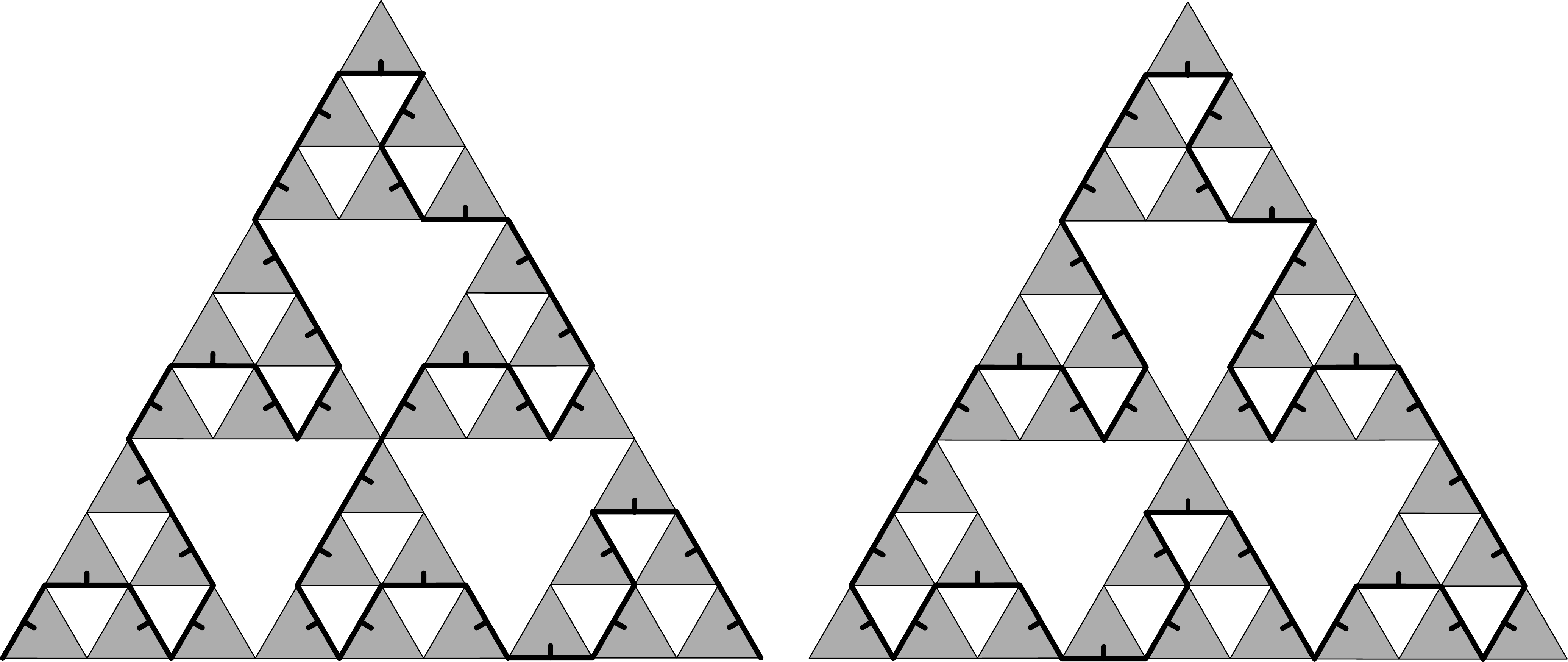}
\centerline{}
\centerline{\bf Figure 3. A tiling-path $(S_{3,2})$ and a tiling-cycle $(D_{3,2})$}
\centerline{\bf on the overall graph $O_3(2)$.}
\end{figure}

\clearpage

\section{The Inscribed Graph $I_n^k$}

Let us define the \emph{Inscribed Graph $(I_n^k)$}, related to the $F_n(k)$ fractal approximating pattern, where we replace all the dark tiles with their centroids as the nodes of the graph and by connecting all the centroids between node-neighbour dark subtriangles we get the edges of the graph. 

It consists of $T_n^{k-1}$ simple $T_n$ sized independent triangular grids which do not share grid points with each other. Connecting edges from the previous approximations remain among the simple triangular grids, otherwise grid points become simple new triangular grids among the connecting edges. See \emph{Figure 4}. 

In this section we will describe the numbers of their nodes and edges, and we will use \emph{Z-paths} and \emph{D-cycles} to calculate the numbers of their Hamiltonian-paths and Hamiltonian-cycles.

\begin{figure}[ht]
\centering
\includegraphics[width=0.6\textwidth]{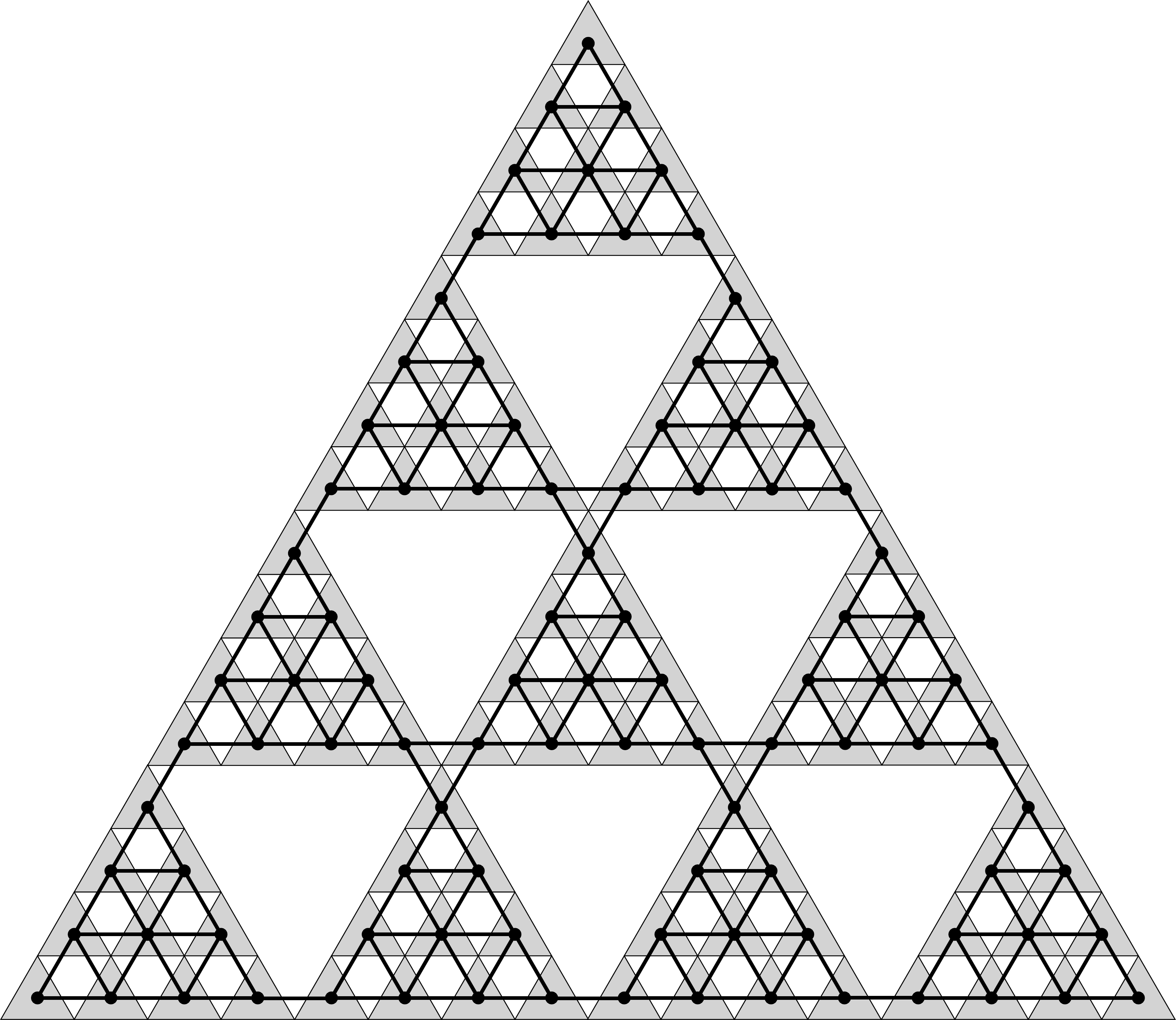}
\centerline{}
\centerline{\bf Figure 4. The Inscribed Graph $I_4^2$}
\end{figure}

\subsection{The numbers of the nodes and the edges in $I_n^k$}

This structure consists of $T_n^k$ nodes in the $k$-th approximation, the $k$-th power of a triangular number: $\left\vert I_n^k \right\vert=T_n^k$.

The number of the edges: $$E(I_n^k)=\sum_{i=1}^{k}3T_{n-1} \cdot T_n^{i-1}$$

because connecting edges from the previous approximations remain among the simple triangular grids.

\subsection{Hamiltonian-paths on $I_n^k$}

We denote the Hamiltonian-paths on the Inscribed Grid by $H_{n,k}$. First we observe special cases, then we find the general formula  to calculate their cardinality.

On the inscribed grid the number of the Hamiltonian-paths is $H_{n,1}=H_n$.

\subsubsection{Hanoi Graph $H_2,k=1$}

$I_2^k$ with all its possible connecting edges is known as the \emph{Hanoi Graph} [H86]. It has a unique Hamiltonian-path (from the leftmost to the rightmost grid point) and a unique Hamiltonian-cycle in any $k$-th approximation which shows how to solve \emph{Hanoi Tower} puzzle if we have $n+1$ pegs in one row and $k$ discs, and only one disc can be moved at one time to a neighbour peg. Discs can be located only in descending order of the disc-sizes.

$H_2=W_2=1$, therefore this is the unique recursive curve on the inscribed graph, the node-rewriting \emph{Sierpi\'{n}ski Arrowhead Curve}, which is also the unique symmetric one.

Hanoi Graphs in general case with more than 3 pegs are also known, but they have other structures than our $I_n^k$ graphs for $n>2$. Our $I_n^k$ graphs can also be represented as all the numbers in $k$ places in the base $T_n$ numeral system.

\subsubsection{Paths on $I_n^2$ and the v-shaped connecting edges}

As we observe the structure of $I_n^k$ for $k>1$ we can see that we have to use W-paths to connect the smaller neighbouring grids to each other. On the sides of the inscribed graph we always have only one possible connecting edge among the smaller grids.

Inside the graph we always have three possible connecting edges between three neighbour grids which form a little triangle, facing downwards. The whole structure looks like a combination of paths and tiles. We have only one entering and one exiting point on the smaller grids, therefore we have to follow their order as the connecting W-path leads the edges among them.

There are no more passages between the smaller grids, but inside the graph, connecting edges of the W-path can take over a corner point from a neighbouring smaller grid. These corner points change the connecting edge to an edge-pair, forming a little v-shaped connection. It modifies the H-path of the smaller grid to a Z-path. The permutation of the smaller grids stands, and the connecting edges follow the W-path with this little modification. By forgetting v-shapes we calculate the number of the possible covering paths on a smaller grid to get $H_{n,2}$ values, denoted by $Y_n$, where $Y_n=H_n+Z_n$.

\subsubsection{Calculating the Hamiltonian-paths on $H_{3,k}$}

We don't have connecting v-shapes and Z-paths on $I_2^k$. They  appear first when $n=3$, and in this case it is easy to calculate them in any order of $k$, because on $I_3^2$ graph our two possible connecting W-paths can always modify exactly one edge to a v-shape, unlike other W-paths in larger orders of $n$.

\smallskip

$H_{3,k}=W_3^a \cdot H_3^b \cdot Y_3^c=2^a \cdot 2^b \cdot 5^c$ (by \emph{Table 1}) where $k>1$, $a=\sum_{i=0}^{k-2}T_3^i$, $b=T_3^{k-1}-a$, $c=a$, therefore we can simplify our formula to

\begin{center}

$H_{3,k}=10^a \cdot 2^{6^{k-1}-a}$ \qquad where \qquad $a=\sum_{i=0}^{k-2}6^i$

\bigskip

$(H_{3,2}=10 \cdot 2^5, \qquad H_{3,3}=10^7 \cdot 2^{29}, \qquad H_{3,4}=10^{43} \cdot 2^{173}, \qquad ...)$

\end{center}

\subsubsection{Calculating the Hamiltonian-paths on $H_{n,2}$}

We find another problem related to v-shaped connecting edges when $n>3$. W-paths have different properties in the same order of $n$. For example, on $I_3^2$ graph our two possible W-paths can always modify exactly one edge to a v-shape.

On $I_5^2$ we have 16 possible W-paths, of which 2 can modify 4 inner edges, 8 can modify 5 inner edges and 6 can modify 6 inner edges to a v-shape, therefore it is difficult to calculate the number of the Hamiltonian paths ($H_{n,k}$) when $n>3$. 

See \emph{Table 4}, where W-paths $(W_n)$ for each $n$ can be separated by the possible number of the connecting v-shapes ($c_m$) to $m$ groups ($g_m$), where  $W_n=\sum_{m=1}^{n-2}g_m$.

\renewcommand{\arraystretch}{1.1}
\[
\begin{tabular}{|c||c|c|c|c|c|}
\hline
\bf
{\small $n$} & \bf { \small $T_n$} & \bf {\small $g_m$} & \bf {\small $b_m$ } & {\small $c_m$ } & {\small $W_n$ } \\ \hline \hline
\bf 3 & 6 & 2 & 5 & 1 & 2 \\ \hline
\bf 4 & 10 & 2 & 8 & 2 & 4  \\ 
\bf & & 2 & 7 & 3 &   \\ \hline
\bf 5 & 15 & 2 & 11 & 4 & 16  \\ 
\bf & & 8 & 10 & 5 &   \\ 
\bf & & 6 & 9 & 6 &   \\ \hline

\bf 6 & 21 & 4 & 14 & 7 & 68  \\ 
\bf & & 22 & 13 & 8 &   \\ 
\bf & & 32 & 12 & 9 &   \\ 
\bf & & 10 & 11 & 10 &   \\ \hline

\bf 7 & 28 & 8 & 27 & 11 & 464  \\ 
\bf & & 76 & 26 & 12 &   \\ 
\bf & & 180 & 25 & 13 &   \\ 
\bf & & 160 & 24 & 14 &   \\ 
\bf & & 40 & 23 & 15 &   \\ \hline

\hline
\end{tabular}
\]

\smallskip

\centerline{\bf \emph{Table 4.} The number of the v-shapes ($c_m$) for $W_n$ values}

\bigskip 
\centerline{The number of the Hamiltonian-paths on $I_n^2$ (for $n>2$) is:}

$$H_{n,2}=\sum_{m=1}^{n-2}g_m \cdot H_n^{b_m} \cdot Y_n^{c_m}$$

$T_n$ little triangular grid can be covered by $b_m$ H-paths and $c_m$ triangular or trapezoid covering paths ($Y_n=H_n+Z_n$). $T_n=b_m + c_m$ for all $m$. 

\begin{remark}
We have left the order number $n$ from $g_m$, $b_m$ and $c_m$ for the sake of simplicity. Self-evidently they always belong to the actual $W_n$ values.
\end{remark}

See \emph{Figure 5} for a possible Hamiltonian-path on the inscribed graph $I_4^2$ and for the explicit formula of $H_{4,2}$.

\begin{figure}[ht]
\centering
\includegraphics[width=0.6\textwidth]{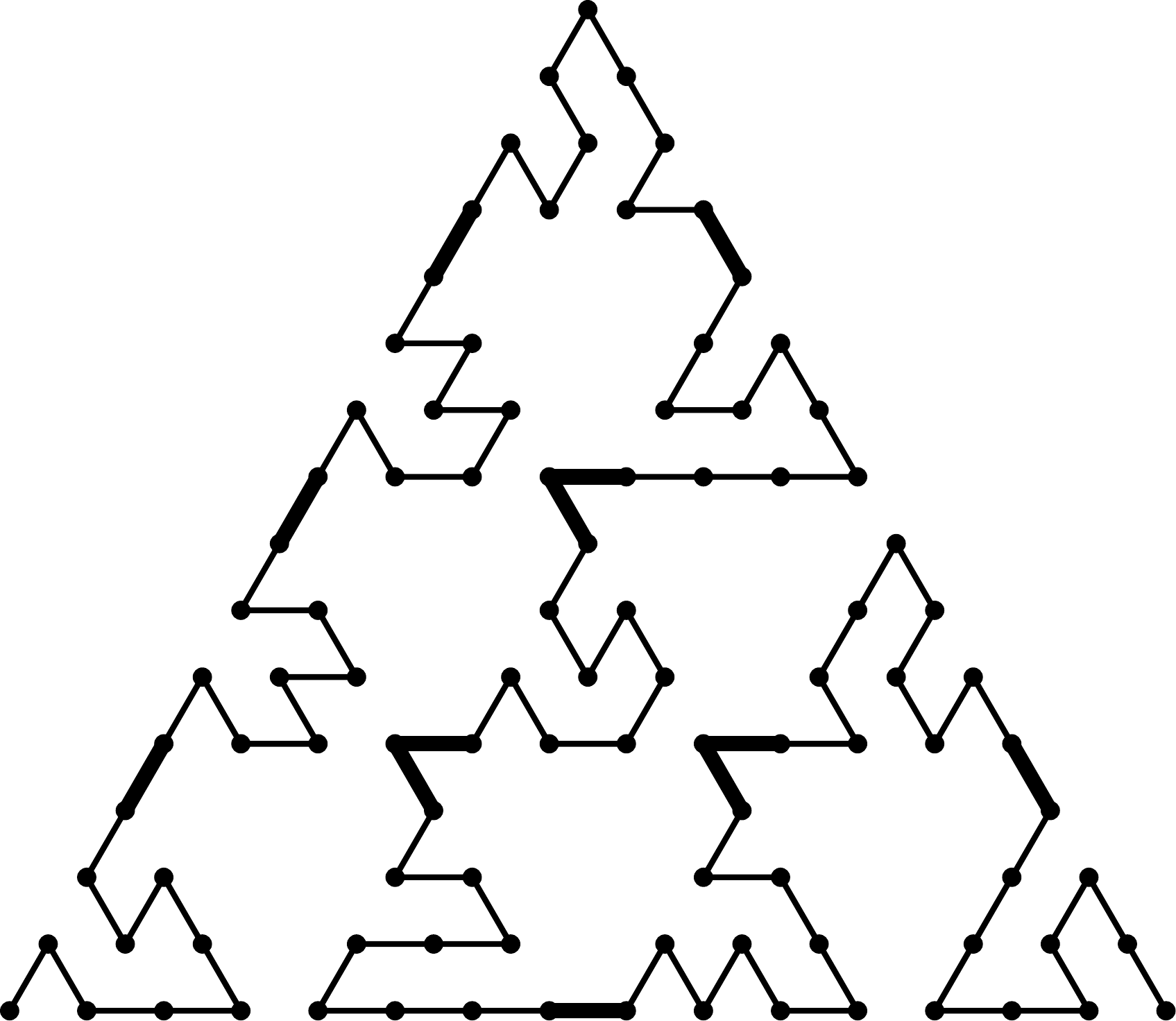}

\bigskip

\centerline{\bf Figure 5. A Hamiltonian-path with v-shaped}
\centerline{\bf connecting edge-pairs and Z-paths on $I_4^2$ graph.}
\smallskip
\centerline{\bf $H_{4,2}=g_1 \cdot H_4^{b_1} \cdot Y_4^{c_1} + g_2 \cdot H_4^{b_2} \cdot Y_4^{c_2}=$}
\smallskip
\centerline{\bf $2 \cdot 10^8 \cdot 21^2 + 2 \cdot 10^7 \cdot 21^3 =  273420000000$}
\end{figure}

\subsubsection{General formula for any $H_{n,k}$}

In the previous subsections we have given the explicit formulas to calculate $H_{n,1}$, $H_{n,2}$, $H_{2,k}$ and $H_{3,k}$ values. Now we give a recursive formula to calculate the Hamiltonian-paths ($H_{n,k}$) on larger Inscribed Graphs $I_n^k$, where $n>3$, $k>2$.

Consider a second order approximation as we can see in \emph{Figure 5}. In the third approximation we can use this second order path like a larger triangular tile, or as in Z-paths by erasing the uppermost grid point as a trapezoid tile, which means we can cover its uppermost first order triangular grid with an H-path or also a Z-path. It modifies our last result (one more first order grid covered by a $Y_n$ path instead of a first order grid covered by an H-path).

We give a recursive formula to calculate the Hamiltonian-paths in general case on $I_n^k$:

\smallskip
$$H_{n,k}= \sum_{m=1}^{n-2}g_m \cdot H_{n,k-1}^{b_m} \cdot \bigg( \dfrac{H_{n,k-1} \cdot Y_n}{H_n} \bigg)^{c_m} $$

\subsection{Hamiltonian-cycles on $I_n^k$}

Here we observe the most complicated structure in this paper, and we give a formula to calculate the number of the Hamiltonian-cycles $(C_{n,k})$ on the Inscribed Graph.

The number of the Hamiltonian-cycles, if $k=1$ is $C_{n,1}=C_n$. For $C_n$ values see \emph{Table 1}, [P14] or [OEIS1].

For the second approximation $(k=2)$ we have to use well-formed connections between $T_n$ smaller triangular grids, therefore we have to use D-cycles. Like in the previous case, inner connecting edges can be substituted by little v-shaped edge-pairs, but their numbers are different for a constant $n$, also for D-cycles.

Consider $n=5$. For $k=1$, $C_{5,1}=C_5=26$. For $k=2$, $D_5=8$, but 6 of the cycles use 4 inner connecting edges, 2 of the cycles use 6 inner connecting edges. Generally on an $I_n^2$ inscribed graph we have $T_n$ smaller grids. The number of their possible connections is $D_n$, which can be separated into $m$ groups ($f_m$), where $D_n=\sum_{m=1}^{n-3}f_m$ and $T_n=r_m+t_m$.

See \emph{Table 5} for $D_n$ values grouped by the number of the v-shapes ($t_m$).

\[
\begin{tabular}{|c||c|c|c|c|c|}
\hline
\bf
{\small $n$} & \bf { \small $T_n$} & \bf {\small $f_m$} & \bf {\small $r_m$ } & {\small $t_m$ } & {\small $D_n$ } \\ \hline \hline
\bf 4 & 10 & 3 & 8 & 2 & 3 \\ \hline
\bf 5 & 15 & 6 & 11 & 4 & 8  \\ 
\bf & & 2 & 9 & 6 &   \\ \hline

\bf 6 & 21 & 6 & 14 & 7 & 42  \\ 
\bf & & 30 & 13 & 8 &   \\ 
\bf & & 6 & 11 & 10 &   \\ \hline

\bf 7 & 28 & 24 & 17 & 11 & 240  \\ 
\bf & & 108 & 16 & 12 &   \\ 
\bf & & 24 & 15 & 13 &   \\ 
\bf & & 84 & 14 & 14 &   \\ \hline

\bf 8 & 36 & 72 & 20 & 16 & 2120  \\ 
\bf & & 432 & 19 & 17 &   \\ 
\bf & & 932 & 18 & 18 &   \\ 
\bf & & 240 & 17 & 19 &   \\ 
\bf & & 444 & 16 & 20 &   \\ \hline

\hline
\end{tabular}
\]

\smallskip
\centerline{\bf \emph{Table 5.} The number of the v-shapes ($t_m$) for $D_n$ values}

\smallskip

The Hanoi Graph ($I_2^k$) has only one Hamiltonian-cycle in any approximations, therefore $C_2,k=1$.

$C_{3,k}=H_{3,k-1}^6$ because $D_3=1$, otherwise there is only one way to connect all the small grids to a cycle, and these cycles do not contain v-shaped edges. The center of the graph is never connected. It is the same as connecting six Hamiltonian-paths ($H_{3,k-1}$) from the previous approximation with each other.

Hamiltonian-cycles on the second order approximations $(C_{n,2})$ can be calculated on the following way, where $n>3$:

$$C_{n,2}=\sum_{m=1}^{n-3}f_m \cdot H_n^{r_m} \cdot Y_n^{t_m}$$ 

\begin{remark}
We have left the order number $n$ from $f_m$, $r_m$ and $t_m$ for the sake of simplicity. Self-evidently they always belong to the actual $D_n$ values.
\end{remark}

See \emph{Figure 6} for a Hamiltonian-cycle on $I_4^2$. 

\begin{figure}[ht]
\centering
\includegraphics[width=0.6\textwidth]{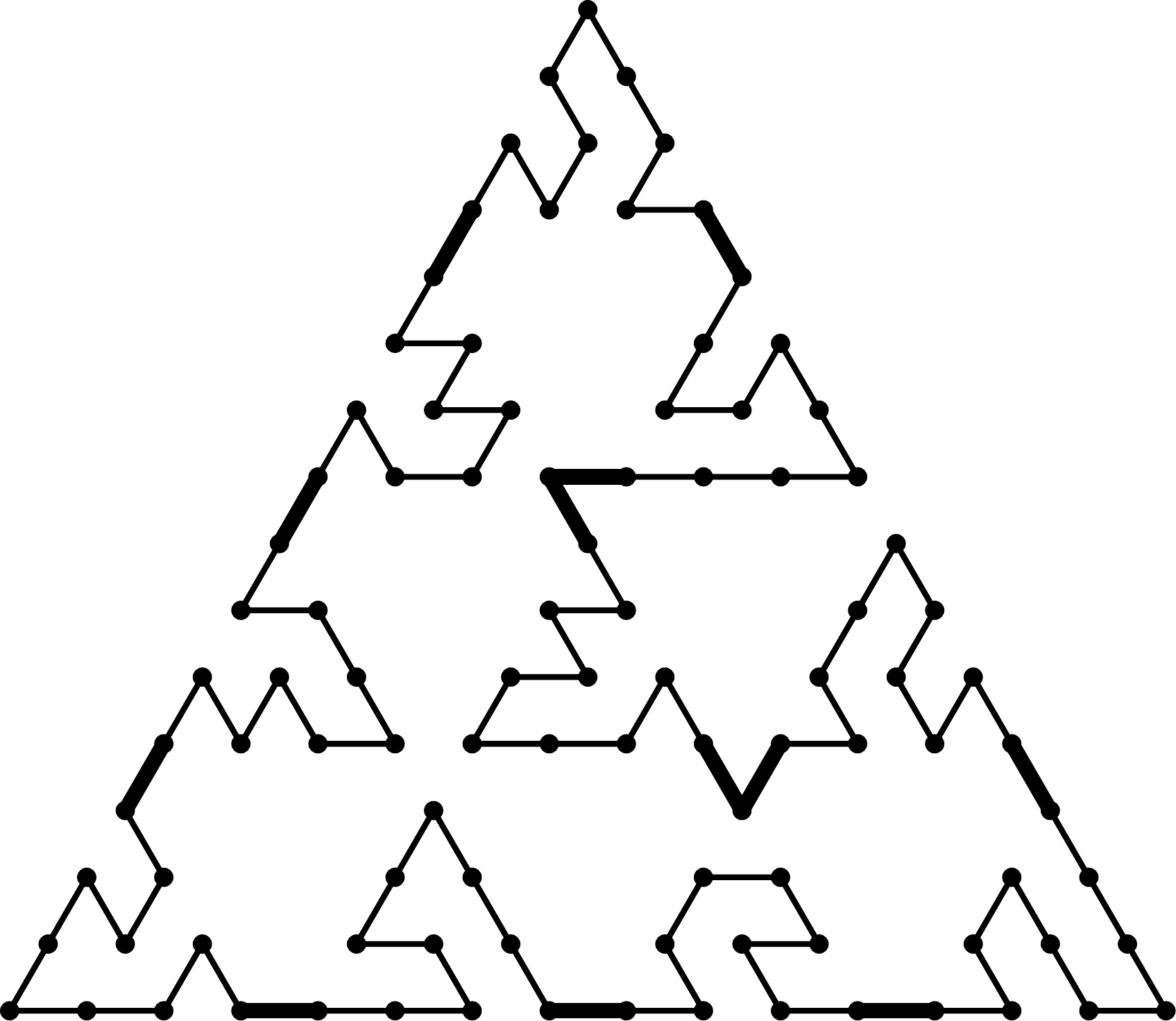}

\bigskip

\centerline{\bf Figure 6. Example of a Hamiltonian-cycle $C_{4,2}$}
\smallskip
\centerline{\bf with v-shaped connecting edge-pairs and Z-paths on $I_4^2$ graph.}
\end{figure}

The recursive formula to calculate the Hamiltonian-cycles in general case on $I_n^k$ is the following:

$$C_{n,k}= \sum_{m=1}^{n-3}f_m \cdot H_{n,k-1}^{r_m} \cdot \bigg( \dfrac{H_{n,k-1} \cdot Y_n}{H_n} \bigg)^{t_m} $$

\subsection{Well-formed Hamiltonian-paths and -cycles}

Hamiltonian-paths and Hamiltonian-cycles on $I_n^k$ have subsets consisting of only well-formed turns (no v-shapes, no Z-paths and H-paths). These paths and cycles on the Inscribed Graph $(I_n^k)$ have the same cardinality as tiling-paths $(S_{n,k})$ and tiling-cycles $(D_{n,k})$ on the Overall Graph $(O_n(k))$.

\section* {Summary}

In \emph{Section 1} we have summarized and supplemented our previous results, which related to the triangular grids, the \emph{generalized Sierpi\'{n}ski Gasket} and the \emph{generalized Sierpi\'{n}ski Arrowhead Curve}. In \emph{Sections 2} and \emph{3} we have observed the background structures of the same fractal family as two kinds of fractal approximating graphs. We have given explicit and recursive formulas to calculate the cardinality of their nodes and edges in both cases, their tiling-paths and -cycles, which cover all the dark tiles of the \emph{Overall Graph}, and their edge covering paths and cycles on the \emph{Inscribed Graph}. We have found their interesting properties, and we have also found new integer sequences.

Our earlier papers also complete this field with some details \emph{[HK15,HK16,K17a]}. For all of my papers please check my Google Scholar site \emph{[KA]}.

\section* {References}

\smallskip \noindent [K17a] Kaszanyitzky, A.: \emph{The generalized Sierpi\'{n}ski Arrowhead Curve} \\
\url{https://arxiv.org/abs/1710.08480}, 2017.
\smallskip

\smallskip \noindent [SEH05] Staji\'c, J., Elezovi\'c-Had\v{z}i\'c, S.: \emph{Hamiltonian walks on Sierpinski and n-simplex fractals}, \url{https://arxiv.org/abs/cond-mat/0310777}, 2005.
\smallskip 

\smallskip \noindent [P14] Pettersson, V. H.: \emph{Enumerating Hamiltonian Cycles}, 2014.
\smallskip

\smallskip \noindent [OEIS1] Sloane, N.J.A.: \emph{The On-line Encyclopedia of Integer Sequences}, 
\\
\url{https://oeis.org/A112676}
\smallskip

\smallskip \noindent [OEIS2] Sloane, N.J.A.: \emph{The On-line Encyclopedia of Integer Sequences}, \\
\url{https://oeis.org/A067771}
\smallskip

\smallskip \noindent [CC06] Chang, S., Chen, L.: \emph{Spanning trees on the Sierpinski gasket} \\
\url{https://arxiv.org/abs/cond-mat/0609453}, 2006.
\smallskip

\smallskip \noindent [H86] \emph{Wolfram MathWorld / Hanoi Graph},
\\ \url{http://mathworld.wolfram.com/HanoiGraph.html}, 1986.
\smallskip

\smallskip \noindent [HK15] Hujter, M., Kaszanyitzky, A.: \emph{Hamiltonian paths on directed grids}, \\ \url{https://arxiv.org/abs/1512.00718}, 2015.
\smallskip

\smallskip \noindent [HK16] Hujter, M., Kaszanyitzky, A.: \emph{Symmetries related to domino tilings on a chessboard}, \url{https://arxiv.org/abs/1603.03298v2}, or in: \it Symmetry: Culture and Science, \sl Vol. 27, pp. 3-10, 2016.
\smallskip

\smallskip \noindent \rm [KA] \emph{Google Scholar citations of Kaszanyitzky, A.},
\\ \url{https://scholar.google.hu/citations?user=i5daxSoAAAAJ}
\smallskip

\end{document}